\documentclass[11pt]{article}
\usepackage{latexsym,amsbsy,amssymb}
\usepackage[latin1]{inputenc}
\usepackage{a4wide}

\newcommand{\R}{\mathbb R}

\newcommand{\g}{\ensuremath{\hat{\mathcal G}} }
\newcommand{\e}{\ensuremath{\hat{\mathcal E}} }
\newcommand{\m}{\ensuremath{\hat{\mathcal E}_M} }
\newcommand{\n}{\ensuremath{\hat{\mathcal N}} }
\newcommand{\Om}{\Omega}
\newcommand{\vphi}{\varphi}\newcommand{\al}{\alpha}
\newcommand{\D}{{\mathcal D}}\newcommand{\pa}{\partial}
\newcommand{\Fl}{\mathrm{Fl}}

\begin{document}
\title{A geometric approach to full Colombeau algebras}
\author{R. Steinbauer\footnote{This work was supported by projects P-16742 and Y-237 of the Austrian Science Fund.}\\
Faculty of Mathematics, University of Vienna\\ 
 Nordbergstr.\ 15, A-1090 Wien, Austria\\
 E-mail: roland.steinbauer@univie.ac.at}
\date{October 10, 2007}
\maketitle
\abstract{We present a geometric approach to diffeomorphism invariant full Colombeau algebras 
which allows a particularly clear view on the construction of the intrinsically defined algebra 
$\g(M)$ on the manifold $M$ given in \cite{gksv}.


}
\section{Introduction.}\label{s1}
In the early 1980-ies J. F. Colombeau \cite{c1}, \cite{c2} constructed algebras of generalised functions
con\-taining the vector space $\D'$ of distributions as a subspace and the space of ${\mathcal C}^\infty$-functions as a 
subalgebra. As associative and commutative differential algebras they combine a maximum of favourable differential 
algebraic properties with a maximum of consistency properties with respect to classical analysis, according to the 
Schwartz impossibility result (\cite{s}). Colombeau algebras since then have proved to be a useful tool in nonlinear 
analysis, in particular in nonlinear PDE with non-smooth data (see \cite{o} and references therein) and have 
increasingly been used in geometric applications such as Lie group analysis of differential equations (see, e.g.\
\cite{dkp}) and general relativity (for a recent review see \cite{sv}).

In this work we shall be exclusively interested in so called {\em full} Colombeau algebras which possess a
canonical embedding of distributions. One drawback of the early constructions (given e.g.\ in \cite{c1}) 
was their lack of diffeomorphism invariance. (Note that this is in contrast to the situation
for the so-called special algebras, which do not possess a canonical embedding of $\D'$: The definition
of special algebras is automatically diffeomorphism invariant, hence these algebras lend themselves 
naturally to geometric constructions, see e.g.\ \cite[Ch.\ 3.2]{gkos}. For matters of embedding $\D'$ into 
simplified algebras in the manifold context see \cite[Ch.\ 3.2.2]{gkos}.) 
It was only after a considerable effort that the full construction could be suitably modified to obtain diffeomorphism invariance (\cite{cm,jel,gfks}). Moreover, the construction of a (full) Colombeau algebra 
$\g(M)$ on a manifold $M$ exclusively using intrinsically defined building blocks was given in \cite{gksv}. 
Note that such an intrinsic construction is vital for applications in a geometric context, e.g.\ in relativity.
The natural next step, i.e., the construction of generalised tensor fields on manifolds again proved to be 
rather challenging. For details we refer to the forthcoming work \cite{gksv2}.
Among other tasks (that are in part discussed in \cite{g}) it was necessary to develop a new 
point of view on the construction given in \cite{gksv}, in particular, on the property that the Lie 
derivative commutes with the embedding. 

In this short note we elaborate on this new point of view by presenting a geometric approach to the 
construction of the algebra $\g(M)$ of \cite{gksv}. We believe that this novel approach serves two purposes: 
It provides a short introduction into diffeomorphism invariant full Colombeau algebras for 
readers not familiar with this topic and it suggests to the experts a very useful shift of focus.
\bigskip

To set the stage for our main topic to be presented in section \ref{s2} we begin by recalling  
the general characteristics of Colombeau's construction on open subsets of $\R^n$. 
It provides associative and commutative differential algebras---from now on denoted by ${\mathcal G}$---satisfying the 
following distinguishing properties:
\begin{enumerate}
 \item[(i)] There exists a linear embedding $\iota:\D'\hookrightarrow{\mathcal G}$ and the function $f(x)=1$ 
  is the unit in ${\mathcal G}$.
\item[(ii)] There exist derivative operators $\pa:\ {\mathcal G}\to{\mathcal G}$ that are linear and satisfy the
 Leibniz rule.
\item[(iii)] The operators $\pa$, when restricted to $\D'$, coincide with the usual partial derivatives, i.e.,
 $\pa\circ\iota=\iota\circ\pa$.
\item[(iv)] Multiplication in ${\mathcal G}$, when restricted to ${\mathcal C}^\infty\times{\mathcal C}^\infty$, 
 coincides with the usual product of functions.
\end{enumerate}
Recall that these properties are optimal in the light of the impossibility result of L. Schwartz (\cite{s}).
Roughly the construction consists of the following steps (for a more elaborate scheme of construction see
\cite[Ch.\ 3]{gfks}):
\begin{enumerate}
 \item [(A)] Definition of a basic space ${\mathcal E}$ that is an algebra, together with linear embeddings
  $\iota:{\mathcal D}'\hookrightarrow{\mathcal E}$ and $\sigma:{\mathcal C}^\infty\hookrightarrow{\mathcal E}$
  where $\sigma$ is an algebra homomorphism. Definition of derivative operators on ${\mathcal E}$ that coincide 
  with the usual derivatives on ${\mathcal D}'$.
 \item [(B)] Definition of the spaces ${\mathcal E}_M$ of moderate and ${\mathcal N}$ of negligible elements 
  of the basic space ${\mathcal E}$ such that ${\mathcal E}_M$ is a subalgebra of ${\mathcal E}$
  and ${\mathcal N}$ is an ideal in ${\mathcal E}_M$ that contains $(\iota-\sigma)({\mathcal C}^\infty)$.
  Definition of the algebra as the quotient ${\mathcal G}:={\mathcal E}_M/{\mathcal N}$.
\end{enumerate}

Many different versions of this construction have appeared over the years, adapted to special
situations or specific applications. In the next section we shall focus on step (A) above in a geometrical 
context.

\section{A geometric approach to Colombeau algebras.}\label{s2}
In this main part of our work we present a geometric approach to diffeomorphism invariant full 
Colombeau algebras. We will carry out our construction on an oriented paracompact smooth Hausdorff 
manifold $M$ of dimension $n$ and proceed in three steps.
\medskip

\noindent
\subsection{The basic space and the embeddings.}
 We want to embed the space of distributions $\D'(M):=(\Omega^n_c(M))'$ (with $\Om^n_c$ denoting the space
 of compactly supported $n$-forms) as well as ${\mathcal C}^\infty(M)$ 
 into our forthcoming basic space $\e(M)$: a  natural choice therefore would be 
 ${\mathcal C}^\infty(\Omega^n_c(M)\times M)$. However, for technical reasons one actually restricts the 
 first slot to elements of $\Omega^n_c$ with unit integral. Denoting their space by $\hat{\mathcal A}_0(M)$ we
 define the basic space by
 \[\e(M):={\mathcal C}^\infty(\hat{\mathcal A}_0(M)\times M).\]
 Elements of the basic space will be denoted by $R$ and its arguments by $\omega$ and $p$.
 Now it is natural to define the embeddings $\sigma$ and $\iota$ by
 \[ \sigma(f)(\omega,p):=f(p)\quad \mbox{resp.\ }\quad \iota(u)(\omega,p):=\langle u,\omega\rangle.\]
 Note that we clearly have $\sigma(fg)=\sigma(f)\sigma(g)$. On the other hand, the formula for $\iota$ might seem
 a little unusual to those acquainted with the works of Colombeau \cite{c1,c2}. 
 In the local situation---that is, on an open subset $\Omega$ of $\R^n$---it has been used by Jelinek 
 in \cite{jel} while the more familiar formula (found e.g.\ in \cite{c2}, however, with an additional reflection) 
 for embedding a distribution $u$ on $\Omega$ is ($\vphi\in\D(\Om)$, $x\in\Om$)
 \begin{equation}\label{colemb}
   \iota(u)(\vphi,x):=(u*\vphi)(x)=\langle u,\vphi(.-x)\rangle.
 \end{equation}
 In the local situation these two formulae give rise to different formulations
 (``Jelinek formalism'' vs.\ ``Colombeau formalism'') which in turn have been shown to give rise to equivalent
 theories in \cite[Ch.\ 5]{gfks}.
(As a historical remark we mention that, in fact, the embedding of distributions in 
Colombeau's first approach (\cite{c0}), if written out explicitly, would take the form 
$\iota(T)(\omega)=\langle T,\omega\rangle$, thus anticipating ``half'' of Jelinek's 
formula. For a detailed discussion of the algebra introduced in \cite[Def.\ 3.4.6]{c0}, see 
\cite[Ch.\ 1.6]{gkos}.)

 In the ``Colombeau formalism'' property (iii) is an easy consequence of the interplay between convolution and 
 differentiation. Indeed, for any multi-index $\al$ we have 
 $(\iota\pa^\al u)(\vphi,x)=(\pa^\al u* \vphi)(x)=\pa^\al(u * \vphi)(x)=(\pa^\al\iota(u))(\vphi,x)$. 

 However, on a manifold---in the absence of a concept of convolution---no analogue to formula (\ref{colemb}) 
 exists and Jelinek's embedding becomes the only one possible. Also on $M$ partial
 derivatives have to be replaced by Lie derivatives with respect to smooth 
 vector fields, which---to keep the presentation simple---we will generally assume to be complete.
 Recall that the Lie derivative of a smooth function $f$ on $M$ w.r.t.\ the 
 vector field $X$ is defined by 
 $L_Xf:=\frac{d}{d\tau}\mathop{\mid}_{\tau=0}\big(\Fl_\tau^X\big)^*f,$
 where $(\Fl_\tau^X)^*$ denotes the pullback under the flow of $X$. 
 So, in order to define derivative operators on our basic space, we first have to introduce 
 the pullback action of diffeomorphisms induced on $\e(M)$, which we shall do next.
 \medskip

\noindent
\subsection{Action of diffeomorphisms.}
 Let $\mu:M\to M$ denote a diffeomorphism of the manifold $M$. Classically we have the following action
 of $\mu$ on smooth functions resp.\ distributions
 \[\mu^*f(p):=f(\mu p)\quad \mbox{and}\quad \langle\mu^*u,\omega\rangle:=\langle u,\mu_*\omega\rangle,\]
 where we have used $\mu p$ as shorthand notation for $\mu(p)$ and $\mu_*\omega$ denotes the pushforward of 
 the $n$-form $\omega$: when written in coordinates
 the second of the above formulae takes the familiar form of the usual pullback of distributions. 
 Therefore the natural definition of the action of $\mu$ on $\e(M)$ which we denote
 by $\hat\mu^*$ is given by
 \[\hat\mu^*R(\omega,p):=R(\mu_*\omega,\mu p).\]
 This definition guarantees that the embeddings are diffeomorphism invariant. Indeed we have
 \[ \hat\mu^*\circ\sigma=\sigma\circ\mu^*\quad \mbox{and}\quad \hat\mu^*\circ\iota=\iota\circ\mu^*,\]
 by the following simple calculation 
 \[
  \hat\mu^*(\iota(u))(\omega,p)=\iota(u)(\mu_*\omega,\mu p)=\langle u,\mu_*\omega\rangle=
  \langle\mu^*u,\omega\rangle=\iota(\mu^*(u))(\omega,p),
 \]
 and similarly for $\sigma$.
 We are now ready to define our derivative operators on $\e(M)$.
 \medskip

\noindent
\subsection{Lie derivatives.}
 Let $X$ be a smooth vector field on $M$ and $R\in\e(M)$. We define the Lie derivative of $R$ w.r.t.\ $X$ by
 \[\hat L_X R:=\left.\frac{d}{d\tau}\right|_{\tau=0}\big(\widehat{\Fl^X_\tau})^*\, R.\]
 Note that property (iii) is now immediate. In fact, the Lie derivative automatically commutes with the 
 embeddings since the action of a diffeomorphism does, i.e., we have
 \[\hat L_X\circ\sigma=\sigma\circ L_X\quad \mbox{and}\quad \hat L_X\circ\iota=\iota\circ L_X.\]
 Rather than giving the one-line proof---which would actually be little more than a replication of the
 above calculation showing that $\iota$ commutes with $\mu^*$---we derive an explicit formula for $\hat L_X R$.
 We have
 \begin{eqnarray}\label{jelder}
  \hat L_XR(\omega,p)
   &=&\left.\frac{d}{d\tau}\right|_{\tau=0}\big(\widehat{\Fl^X_\tau})^*\, R(\omega,p)
   \,=\,\left.\frac{d}{d\tau}\right|_{\tau=0}R\big((\Fl^X_\tau)_*\omega,\Fl^X_\tau p\big)\nonumber\\
   &=&d_1R(\omega,p)\left.\frac{d}{d\tau}\right|_{\tau=0}\underbrace{\big(\Fl^X_\tau\big)_*\omega}_                 {(\Fl^X_{-\tau})^*\omega}
   \,+\,d_2R(\omega,p)\underbrace{\left.\frac{d}{d\tau}\right|_{\tau=0}\Fl^X_\tau p}_{X(p)}\nonumber\\
   &=&-d_1R(\omega,p)\,L_X\omega+L_XR(\omega,.)\mid_p,                  
 \end{eqnarray}
 where in the second line we have used the chain rule. Note that already in the definition of
 $\e(M)$ we have used calculus in infinite dimensions, but it is at this point where we can
 emphasise the inevitability of using it:
 In a manifold setting, there is only ``Jelinek's formalism'' available. To have the embedding 
 commute with the action of diffeomorphisms respectively Lie derivatives, the structure of the 
 very formula (2) clearly necessitates taking derivatives also w.r.t.\ the $\omega$-slot. 
 This is the ultimate reason for requiring smoothness w.r.t.\ all variables in the definition 
 of the basic space. For a detailed account on this matters see \cite[pp.\ 103--107]{gkos}.
 However, in order to keep our presentation
 simple, we only remark that Chapters 4, 6 and 14 of \cite{gfks} provide all relevant
 details as to calculus in convenient vector spaces (see \cite{km}) for the setting at hand, 
 and (cf.\ \cite{cm}, p.\ 362) `we invite the reader to admit the [respective] smoothness properties'. 

 Rather we would like to draw the attention of our readers to the fact that formula (\ref{jelder})---in 
 the local context---already appeared in Remark 22 of \cite{jel}. However, in that reference 
 it is an operational consequence of Jelinek's formalism, whereas here it arises
 as a simple consequence of our natural choice of definitions. 
 \bigskip
 
 At the end of this section we digress to remark that the geometric
 approach to the definition of the Lie derivative 
 can also be employed to define the usual derivative for distributions, say on the real line. 
 Indeed, taking $X=\pa_x$, the flow is given by a translation, i.e., $\Fl^X_\tau x=x+\tau$ and we have
 \begin{eqnarray*}
  \langle u',\vphi\rangle
  &=&\langle L_Xu,\vphi\rangle
  \,=\,\langle\left.\frac{d}{d\tau}\right|_{\tau=0}\big(\Fl^X_\tau)^*u,\vphi\rangle\\
  &=&\left.\frac{d}{d\tau}\right|_{\tau=0}\langle u,\big(\Fl^X_{-\tau}\big)^*\vphi\rangle
  \,=\,\langle u,-L_X\vphi\rangle
  \,=\,-\langle u,\vphi'\rangle.
 \end{eqnarray*}
 Observe the somehow amusing fact that in this approach no explicit reference to integration by parts occurs. 
 Indeed, in this way integration by parts follows from the translation invariance of the integral.

\section{Summary and conclusions.}
 To sum up we have constructed the basic space for the Colombeau algebra 
 $\g(M)$ of \cite{gksv} together with the embeddings of smooth functions resp.\ distributions
 and the Lie derivative. Recall that our definitions were all geometrically motivated 
 and allowed us to obtain properties (i)---(iii) in a remarkably effortless way.

 Finally we sketch how one establishes property (iv), i.e., the constructions outlined in step (B) above. 
 To this end we employ a quotient construction
 to identify the images of ${\mathcal C}^\infty(M)$ under the two embeddings $\sigma$ and $\iota$.
 More precisely, we want to identify the two middle terms of
 \[\sigma(f)(\omega,p)=f(p)\ \sim\ \int f(q)\omega(q)=\iota(f)(\omega,p).\]
 The obvious idea is to set $\omega(q)=\delta_p(q)$ which clearly is only possible asymptotically.
 At this point the usual asymptotic estimates (`tests' in the language of \cite{gfks}) 
 come into play (for the present case see 
 \cite[Defs.\ 3.10, 3.11]{gksv}) and are used to single out the moderate resp.\ the negligible elements of 
 $\e(M)$. Denoting the respective spaces by $\m(M)$ resp.\ $\n(M)$ we finally may define 
 the Colombeau algebra on the manifold $M$ by 
 \[\g(M):=\m(M)/\n(M).\]
 It is a differential algebra with respect to the Lie derivative w.r.t.\ smooth vector fields and
 the Lie derivative commutes with the embedding of distributions. Moreover it localises appropriately
 to the diffeomorphism invariant local algebra $\g^d(\Omega)$ of \cite{gfks}.
 \medskip 

 The new aspect we have demonstrated in this contribution is the following: The fact that the Lie 
 derivative commutes with both embeddings (on the level of the basic space, hence with the embedding
 of distributions after taking the quotient) is a consequence of diffeomorphism invariance of the 
 embeddings which itself is due to a natural choice of the definition of the action of diffeomorphisms 
 on the basic space as well as of the natural definition of the embeddings themselves. Hence in the global setting 
 property (iii)---which in the local context follows from the properties of the convolution---is a direct 
 consequence of the diffeomorphism invariance of the construction.
\bigskip

{\em Acknowledgements:} The author wishes to express his gratitude to the organisers of the 
GF07 conference in B\c{e}dlewo, in particular to Swiet\l{}ana Minczewa-Kami\'nska and to 
Andrzej Kami\'nski.

\end{document}